\documentclass[12pt]{article}

\usepackage{amsmath}
\usepackage{multirow} 
\usepackage{latexsym,amsfonts}
\usepackage{subcaption}
\usepackage{pdflscape}
\usepackage{caption}
\usepackage{enumitem} 
\usepackage{authblk}
\newtheorem{thm}{Theorem}[section]

\newtheorem{prop}[thm]{Proposition}

\newtheorem{Rem}[thm]{Remark}

\usepackage{multicol}

\author[1]{Sanja Rukavina}  
\author[2]{Vladimir D. Tonchev}
\affil[1]{Faculty of Mathematics, University of Rijeka, 51000 Rijeka, Croatia\\

E-mail: sanjar@math.uniri.hr, ORCID: 0000-0003-3365-7925}
\affil[2]{ Department of Mathematical Sciences, Michigan Technological University, Houghton, MI 49931, USA\\

E-mail: tonchev@mtu.edu, ORCID: 0000-0003-1806-3571}
\title{New examples of self-dual near-extremal ternary codes of length 48 derived from 2-(47,23,11) designs}

\date{}
\begin{document} 
\maketitle

\begin{abstract} 
 In a recent paper \cite{AR-HAR}, 
Araya and Harada gave   examples of self-dual near-extremal ternary codes of length 48 for $145$ distinct values 
of the number $A_{12}$ of codewords  of minimum weight 12, and raised the question about the existence of codes for
other values of $A_{12}$.
 In this note, we use symmetric 2-$(47,23,11)$ designs with an automorphism group
 of order 6 to construct
self-dual near-extremal ternary codes of length 48 for $150$ new values of $A_{12}$.
\end{abstract}

\vspace*{0.5cm}



\section{Introduction}

We assume familiarity with the basic facts and notions 
from error-correcting codes and combinatorial designs 
\cite{AK, BJL, HP, Ton21, ton88}.
All codes considered in this paper are ternary.  

The minimum  weight $d$
of a ternary self-dual code of length $n$ divisible by 12
satisfies the upper bound $d\le n/4 + 3$ \cite[9.3]{HP}. A self-dual ternary code of length $n$ divisible by 12 with minimum weight $d$ is  called {\it extremal} if $d=n/4+3$ \cite{HP}, and {\it near-extremal} if $d=n/4$  \cite{AR-HAR, newAM}. 
Any extremal ternary self-dual code supports combinatorial 5-designs
by the Assmus-Mattson theorem \cite{AM},  \cite[8.4]{HP}. It was recently proved by Miezaki, Munemasa, and Nakasora \cite{newAM}
that the supports of all codewords of weight $w\le 6m-3$ in any ternary near-extremal self-dual code of length $n=12m$ are the blocks
of a combinatorial 1-design. Thus, extremal and near-extremal self-dual codes are interesting from both coding and design theoretical 
point of view.

The classification of extremal ternary self-dual codes of length $n$ divisible by $12$ has been completed  only for the lengths $n=12$ and $n=24$.
Up to equivalence, there is one extremal code of length $12$, being the extended ternary Golay code, and there are
 two  extremal codes of length 24:
the extended quadratic-residue code $QR_{23}^*$ \cite{AM} and the Pless symmetry code $C(11)$ \cite{Pless69}, \cite{Pless72}. 
The only known extremal code of length $n=36$ is the Pless symmetry code $C(17)$ \cite{Pless69}, \cite{Pless72}. 
It was shown recently that the Pless symmetry code $C(17)$ is equivalent to a code spanned by the incidence matrix of of a 
symmetric 2-$(36,15,6)$
design with a trivial full automorphism group \cite{Ton21}, as well as by the incidence matrix of a unique
 $2$-$(36,15,6)$ design that admits an involution \cite{RT36}.
Two extremal ternary self-dual codes of length $48$ are known: the extended quadratic-residue code $QR_{47}^*$ and the Pless symmetry code $C(23)$. Finally, three  extremal codes of length $n=60$ are known:  the extended quadratic-residue code $QR_{59}$, the Pless symmetry code $C(29)$,
and a code found by Nebe and Villar \cite{NV}.\\

The sparsity of extremal codes has spurred some recent interest in near-extremal codes.
Araya and Harada \cite{AR-HAR} 
proved  that the number $A_{3i}$ of codewords of weight $3i$, $m \le i \le 4m$,
in a ternary near-extremal self-dual code of length  $n=12m$,
is divisible by 8.
In addition, the number $A_{12}$ of codewords of minimum weight 12 in a ternary near-extremal self-dual code
of length 48
is $A_{12}=8\beta$ for some $\beta$ in the range
$1\le \beta \le 4324$ \cite[page 1831]{AR-HAR}. 
Araya and Harada gave
  examples of near-extremal ternary codes of length 48 for $145$ distinct values of
 $A_{12}$ (sets $\Gamma_{48,1}$ and $\Gamma_{48,2}$ in \cite{AR-HAR}), and raised the question
for the existence of codes for other values of $A_{12}$
\cite[page 1838, Question 1]{AR-HAR}.\\

In the next section, we use symmetric $2$-$(47,23,11)$ designs admitting an automorphism of order six to find many new examples of self-dual near-extremal ternary codes of length $48$. We constructed examples of self-dual near-extremal ternary codes of length 48 with $150$ distinct values of $A_{12}$ not covered in \cite{AR-HAR}.  
In our computations, in addition to our own computer programs, we used the computer programs of V. \'{C}epuli\'{c} \cite{cep}
for the construction of orbit matrices, and the computer algebra system MAGMA \cite{magma} for  computing the codes
and their weight distributions.

\section{Near-extremal [48,24,12] codes derived from symmetric 2-(47,23,11) designs}
\label{sec2}
The following statement gives a simple construction of ternary self-dual codes of length 48.

\begin{thm}
\label{t1}
Let $M$ be a $47\times 47$ $(0,1)$-incidence matrix of a symmetric 2-$(47,23,11)$ design,
and let $G$ be the $47\times 48$ matrix obtained by adding to $M$ the all-one column.
Then the row space of $G$ over $GF(3)$ is a ternary self-dual code of length 48.
\end{thm}
{\bf Proof}.
Since every row of $G$ has weight $24 \equiv 0 \pmod 3$, and the inner product of every two distinct rows of
$G$ is equal to $12 \equiv 0 \pmod 3$, the ternary code $L$ spanned by the rows of $G$ is self-orthogonal.
Since $23$ is not divisible by 3, and $23-11=12$ is divisible by 3, but not divisible by 9, it follows from
\cite[Theorem 4.6.2, (b)]{AK} that the rank of $M$ over $GF(3)$  is equal to $(47+1)/2=24$, hence the code $L$ is self-dual.
$\Box$\\

We note that the two known extremal ternary self-dual codes of length 48, $QR_{47}^*$ and $C(23)$, are obtainable via
the construction of Theorem \ref{t1} from symmetric $2$-$(47,23,11)$ designs associated with the Paley-Hadamard matrices 
of type I and II and order 48, respectively \cite{Ton21}.\\

 In this paper, we use the method for refinement and indexing of orbit matrices for presumed action of an abelian automorphism group \cite{c-r, cep} to construct 2-$(47,23,11)$ designs invariant under the cyclic group $C_6$ of order 6 with  orbit lengths distribution 
 $(1,2,2,3,3,6,6,6,6,6,6)$. There are $32$ orbit matrices for such an action of an automorphism of order six on a 2-$(47,23,11)$ design. 
 Our computations show that only four orbit matrices (given in the Appendix) yield designs which generate near-extremal ternary codes
via the construction of Theorem \ref{t1}.
The results are summarized in Table \ref{tab_designs}, where $d$  denotes the minimum weight of a code.

\begin{table}[ht] 
	\centering
		\begin{tabular}{|c|c|c|c|c|}
		\hline
			Orbit matrix & OM1&OM2&OM3&OM4\\
			\hline
			Non-isomorphic&$70400$&$24576$ &$63488$&$24576$\\
			designs& & & &\\
			\hline
			\# Codes with $d=12$&$43338$&$11884$ &$22698$&$11884$\\
			\hline
            Inequivalent codes with $d=12$&$1662$ &$1073$ &$1200$ &$1073$\\
			\hline
			\# Distinct $A_{12}$ &$165$ &$152$ &$161$ &$152$\\
			\hline
		\end{tabular}\caption{$2$-$(47,23,11)$ designs  and their codes}\label{tab_designs} 
\end{table}

\begin{table}[ht] \begin{footnotesize} 
	\centering
		\begin{tabular}{|c|c|}
		\hline
			$\Gamma_{OM_1}$& $320, 323, 324, 326, 338, 340, 341, 346, 348, 349, 350, 352, 353, ..., 357,$\\
			& $359, 360, ..., 468, 470, 471, ..., 480, 482, 483, ..., 486, 489, 490, ..., 494,$\\
			&$496, 497, ..., 500, 504, 506, 512, 516, 518, 522, 524, 528, 536, 560 $ \\
			\hline
			$\Gamma_{OM_2}$& $313, 329, 331, 332, 333, 334, 337, 338, 339, 343, 344, ..., 349,
			351, 352, ..., 450,$ \\
			& $452, 453, ..., 459, 461, 462, 464, 466, 467, 468, 470, 472,474, 476, 478, 479,$     \\
			&$480,482, 484,486, 488, 494, 496, 500, 503, 504, 506, 512,524,528, 554, 560$ 
			\\
			\hline
			$\Gamma_{OM_3}$&$320, 323, 324, 326, 338, 340, 341, 346, 348, 349, 350, 353, ..., 357, 359,$\\
			& $360, ..., 464, 466, 467, 468, 470, 471, ..., 480, 482, 483, ..., 486, 488, 489, ...,492, $\\
			&$494, 496, 497, ..., 500, 504, 506, 512, 516, 522, 524, 528, 536, 560 $ \\
			\hline
			$\Gamma_{OM_4}$&$313, 329, 331, 332, 333, 334, 337, 338, 339, 343, 344, ..., 349,
			351, 352, ..., 450,$ \\
			& $452, 453, ..., 459, 461, 462, 464, 466, 467, 468, 470, 472,474, 476, 478, 479,$     \\
			&$480,482, 484,486, 488, 494, 496, 500, 503, 504, 506, 512,524,528, 554, 560$ 
			\\
			\hline
		\end{tabular}\caption{  $\Gamma_{OM_i}$, $i=1,2,3,4$} \label{A12} \end{footnotesize}
\end{table}
The numbers $A_{12}$  of minimum weight codewords in the newly found near-extremal codes are given by
$$A_{12} \in \{8\beta | \beta \in \Gamma_{OM_i} \},$$
where the set $\Gamma_{OM_i}$,  $i \in \{1, 2, 3, 4\}$, is given in Table \ref{A12},
and contains all distinct values of $\beta$ for 
the near-extremal ternary codes obtained from designs with orbit matrix $OM_i$.
A list of 181 2-$(47,23,11)$ designs that generate codes with different number of codewords of minimum weight 12
is available at 
\begin{verbatim}
 https://www.math.uniri.hr/~sanjar/structures/
 \end{verbatim}
All these $181$ designs have the cyclic group $C_6$ as the full automorphism group.  
The data from Table \ref{A12} can be summarized as follows.
 
\begin{prop} \label{prop}
 There is a ternary near-extremal self-dual code of length $48$ with $A_{12} \in \{8\beta | \beta \in \Gamma \},$ where $A_{12}$ is the number of codewords of weight $12$ and $\Gamma= \{
 313, 320, 323, 324, 326, 329, 331, 332, 333, 334, 337, ..., 341, 
343,344,\\ ..., 468, 470, 471, ..., 480, 482, 483, ..., 486, 488, 489, ...,494, 496, 497, ..., 500, 503,\\ 504, 506, 512, 516, 518, 
522, 524, 528, 536, 554, 560   \}$.
\end{prop}

 Araya and Harada  \cite{AR-HAR}  found $86$ codes with distinct values   $A_{12}\equiv 0 \pmod {48}$ (see $\Gamma_{48,1}$ in \cite{AR-HAR}, p. 1831), plus 59 codes with distinct values $A_{12}$ not divisible by {48} 
 \cite[$\Gamma_{48,2}$]{AR-HAR}.
  In our list, there are $31$ values of $A_{12}$ divisible by $48$ that are covered by $\Gamma_{48,1}$ in \cite{AR-HAR}. 
Our examples of near-extremal ternary codes for the remaining $150$ distinct values of $A_{12}$ were not previously known, since $A_{12}/8<282$ for all $A_{12}$ covered by $\Gamma_{48,2}$ in \cite{AR-HAR}.  Therefore, as a result of our construction we found examples of $150$ new self-dual near-extremal ternary codes with values of $A_{12}$ which were not previously known, and correspond to
the values of $\beta$ from the set $\Gamma$ in Proposition \ref{prop} which are not divisible by six.

\begin{Rem}
{\rm
In \cite{c-r47}, fifty-four symmetric  $2$-$ (47, 23, 11)$ designs admitting a faithful action of a Frobenius group of order $55$ were
constructed. We computed the ternary codes of these designs and found that fifteen designs yield near-extremal self-dual ternary codes
with $A_{12}\in\{ 1584, 1680, 2640, 3792 \}$.  Codes with these values of $A_{12}$ were previously found in \cite{AR-HAR}. 
 
 }
\end{Rem}

 \section{Appendix}
 
Orbit matrices for the action of an automorphism of order six on 2-$(47,23,11)$ designs  which generate new near-extremal self-dual ternary codes of length $48$ via Theorem \ref{t1}:\\

\begin{minipage}[b]{0.5\linewidth} 
 
	\centering
		\begin{tabular}{c|c}
		$OM_1$& 1  2  2  3  3  6  6  6  6  6  6\\
		\hline
1&0   2   0   3   0   6   6   6   0   0   0\\
  2&  1   2   2   0   0   6   3   0   3   3   3\\
   2& 0   2   0   0   3   3   3   3   6   3   0\\
    3&1   2   2   3   3   2   2   2   2   2   2\\
    3&0   2   0   3   0   2   2   2   4   4   4\\
   6& 1   1   1   1   1   2   3   4   2   5   2\\
   6& 1   1   1   1   1   2   3   4   4   1   4\\
   6& 1   0   0   2   2   4   3   2   3   3   3\\
   6& 0   0   2   2   1   3   3   3   4   3   2\\
   6& 0   1   1   1   2   2   5   2   2   3   4\\
   6& 0   1   1   1   2   4   1   4   2   3   4
		\end{tabular}  
 \end{minipage}
\begin{minipage}[b]{0.5\linewidth} 
 
	\centering
		\begin{tabular}{c|c}
		$OM_2$& 1  2  2  3  3  6  6  6  6  6  6\\
		\hline
   1 &0   2   0   3   0   6   6   6   0   0   0\\
    2&1   2   2   0   0   6   3   0   3   3   3\\
    2&0   2   0   0   3   3   3   3   6   3   0\\
    3&1   2   2   3   3   2   2   2   2   2   2\\
    3&0   2   0   3   0   2   2   2   4   4   4\\
    6&1   1   1   1   1   3   1   5   3   3   3\\
    6&1   1   1   1   1   1   5   3   3   3   3\\
    6&1   0   0   2   2   4   3   2   3   3   3\\
    6&0   0   2   2   1   3   3   3   4   3   2\\
    6&0   1   1   1   2   3   3   3   3   1   5\\
    6&0   1   1   1   2   3   3   3   1   5   3\\
		\end{tabular}  
 
\end{minipage}

\vspace{0.5cm}

 \begin{minipage}[b]{0.5\linewidth} 
  
	\centering
		\begin{tabular}{c|c}
		$OM_3$& 1  2  2  3  3  6  6  6  6  6  6\\
		\hline
1& 0   2   0   3   0   6   6   6   0   0   0\\
 2&   1   2   2   0   0   3   3   3   6   3   0\\
 2&   0   2   0   0   3   6   3   0   3   3   3\\
 3&   1   2   2   3   3   2   2   2   2   2   2\\
 3&   0   2   0   3   0   2   2   2   4   4   4\\
 6&   1   1   1   1   1   2   5   2   2   3   4\\
 6&   1   1   1   1   1   4   1   4   2   3   4\\
 6&   1   0   0   2   2   3   3   3   4   3   2\\
 6&   0   0   2   2   1   4   3   2   3   3   3\\
 6&   0   1   1   1   2   2   3   4   2   5   2\\
 6&   0   1   1   1   2   2   3   4   4   1   4\\
		\end{tabular}  
\end{minipage}
\begin{minipage}[b]{0.5\linewidth} 
  
	\centering
		\begin{tabular}{c|c}
		$OM_4$& 1  2  2  3  3  6  6  6  6  6  6\\
		\hline
1&0   2   0   3   0   6   6   6   0   0   0\\
2&    1   2   2   0   0   3   3   3   6   3   0\\
2&    0   2   0   0   3   6   3   0   3   3   3\\
3&    1   2   2   3   3   2   2   2   2   2   2\\
3&    0   2   0   3   0   2   2   2   4   4   4\\
6&    1   1   1   1   1   3   3   3   3   1   5\\
6&    1   1   1   1   1   3   3   3   1   5   3\\
6&    1   0   0   2   2   3   3   3   4   3   2\\
6&    0   0   2   2   1   4   3   2   3   3   3\\
6&    0   1   1   1   2   3   1   5   3   3   3\\
6&    0   1   1   1   2   1   5   3   3   3   3\\
		\end{tabular}  
 
\end{minipage}

\end{document}